\documentclass [reqno]{amsart}
\usepackage[english]{babel}
\usepackage{textcomp}
\usepackage{amsmath,amsfonts,amssymb,amsthm,mathrsfs,float}
\usepackage{graphicx}
\usepackage{lmodern}
\usepackage{hyperref}
\usepackage{etoolbox}
\usepackage{mathptmx}
\patchcmd{\abstract}{\scshape\abstractname}{\textbf{\abstractname}}{}{}
\hypersetup{urlcolor=blue, colorlinks=true, citecolor=blue}  
\usepackage{cite}

\renewcommand{\Re}{\operatorname{Re}}

\DeclareGraphicsExtensions{.png,.pdf,.eps,.jpg}

\def\Z{\mathbb{Z}}

\def\R{\mathbb{R}}

\def\C{\mathbb{C}}

\def\Abf{\mathbf{A}}
\def\Tbf{\mathbf{T}}

\def\({\Big(}

\title[{\tiny Approximation and Homogenization of thermoelastic wave model}]{
 \bf Approximation and Homogenisation of thermoelastic wave model}
\author[S. Nafiri]{S. Nafiri$^{\dagger}$}
\date{}
\thanks{$^{\dagger}$ Département de Mathématiques, informatique et géomatique, Ecole Hassania des Travaux Publics, Km 7 Route d'El Jadida Casablanca BP 8108, Morocco. Email address: nafiri@ehtp.ac.ma $\mid$ nafirisalem@gmail.com (S. Nafiri).}

\date{}
\begin{document}
\maketitle
\numberwithin{equation}{section}
\newtheorem{thm}{Theorem}[section]
\newtheorem{lem}{Lemma}[section]
\newtheorem{prop}{Proposition}[section]
\newtheorem{Def}{Definition}[section]
\newtheorem{rmk}{Remark}[section]
\newenvironment{dem}{\noindent\textbf{Proof.~}}{\hfill$\blacksquare$\smallbreak}
\begin{center}
\textit{Dedicated to the memory of Professor Hammadi Bouslous}
\end{center}
\begin{abstract}
This paper deals with the approximation and homogenization of thermoelastic
wave model. First, we study the homogenization problem of a weakly
coupled thermoelastic wave model with rapidly varying coefficients, using
a semigroup approach, two-scale convergence method and some variational
techniques. We show that the limit semigroup can be obtained by using a weak
version of the Trotter Kato convergence Theorem. Secondly, we consider the
approximation of two thermoelastic wave model, one with exponential decay
and the other one with polynomial decay. the numerical experiments indicate
that the two discrete systems show different behavior of the spectra. Moreover,
their discrete energies inherit the same behavior of the continuous ones. Finally
we show numerically how the smoothness of data can impact the rate of decay
of the energy associated the weakly coupled thermoelastic wave model.\\\\
{\bf Key words:}
Homogenization, weakly coupled thermoelastic eqn, periodic coefficients, long-time behavior, two scale convergence, exponential stability, polynomial stability, spectral element.
\\\\
{\textbf{Mathematics Subject Classification:} 35B27, 93C20, 93D20, 73C25, 65M06, 65M60, 65M70.}
\end{abstract}
\section{\bf Introduction}
\label{sect1}
A lot of attention has been given to the equations of thermoelasticity over the past half century with motivations from classical natural sciences and applications to materials engineering (e.g., metals, polymers, ceramics, liquid crystals, composites, etc.) Much progress has been made in understanding the transfer of energy in flexible structures, and we refer to \cite{Carlson72,Chadwick60,GRT92} for more details about these topics. For example, on p. 1163 of \cite{GRT92}, the authors mention: "The thermally induced vibrations that hampered the recently launched Hubble space telescope have highlighted the coupling between mechanical vibration and heat transfer and the need to model and control thermoelastic phenomena in flexible structures." 

The equations of thermoelasticity model the change in the size and shape of a material (pure or composite) as the temperature of that material fluctuates, this is called thermal expansion. Materials that are more elastic will expand or contract more than those that are less elastic. Engineers use their understanding of thermoelasticity to design composite materials that can withstand fluctuations in temperature without breaking.

Composite materials are widely used in engineering construction because of their properties
that, in general, are better than those of their individual constituents. Moreover, they are known to display properties often not exhibited by their constitutive parts, , see \cite{Milton2001}. Often we would like to predict the properties of composites without having to build them; mathematically this is a question of solving PDEs with rapidly oscillating coefficients, therefore Homogenization is a natural tool.

The main goal of this chapter is to investigate homogenization properties of linear thermoelasticity.
It is to highlight that as far as homogenization of linear thermoelasticity is concerned, Francfort \cite{Francfort83} was probably the first to investigate the homogenization problem of the following system
\begin{equation}
\label{thermoe0}
 \left\{\begin{array}{llll}
                      c(\frac{x}{\varepsilon})\dfrac{\partial^2 u_\varepsilon}{\partial t^2}=div(a(\frac{x}{\varepsilon})\nabla u_\varepsilon)-\gamma(\frac{x}{\varepsilon})\frac{\partial \theta_\varepsilon}{\partial x}\quad & x\in\Omega,\;t\geqslant 0\\
                      d(\frac{x}{\varepsilon})\dfrac{\partial\theta_\varepsilon}{\partial t}=div(b(\frac{x}{\varepsilon})\nabla \theta_\varepsilon)
		      +\gamma(\frac{x}{\varepsilon}) \dfrac{\partial^2 u_\varepsilon}{\partial t \partial x}\quad &x\in\Omega,\;t\geqslant 0\\
u_\varepsilon=0=\theta_\varepsilon\quad & x\in\partial\Omega,\; t>0\\ 
u_\varepsilon(x,0)=\tilde{u}_{0}(x),\; \dfrac{\partial u_\varepsilon}{\partial t}(x,0)=\tilde{u}_{1}(x),\; \theta_\varepsilon(x,0)=\tilde{\theta}_{0}(x) \quad &x\in \Omega.
                     \end{array}
                   \right.
\end{equation}
For fixed $\varepsilon > 0$, the long-time behaviour of the solutions of \eqref{thermoe0} have been intensively studied in many works under
various assumptions on the coupling operator and initial conditions (see
for example \cite{KBT1997,LR05,LZ99,MN16,N21}).
Following Bensoussan, Lions and Papanicolaou \cite{BLP78} and Sanchez-Palencia \cite{Sanchez80} and by using semigroup approach, the author showed that, as the period of the coefficients goes to zero, the solution of these equations converges to the solution of a related constant coefficient problem, the homogenized problem. In other words: given a family of differential operators $\mathbf{A}_\varepsilon$, assume that the boundary value problem
\begin{equation}
 \left\{\begin{array}{ll}
                      z'_\varepsilon (t)=\mathbf{A}_\varepsilon z_\varepsilon (t)\quad t>0,\\
z_{\varepsilon}(0)=z,
                     \end{array}
                   \right.
\end{equation}
is well-posed in a Sobolev space $H$ for all 
$\varepsilon>0$, and the solutions $z_{\varepsilon}$ form a bounded subset in $H$, therefore there is at least one weak limit $z_{0}$ in $H$ for the family $z_\varepsilon$, when the period $\varepsilon\to 0$, where $z_{0}$ satisfies the homogenized evolution equation
\begin{equation}
\label{eq1.1}
 \left\{\begin{array}{ll}
z'_{0} (t)=\mathbf{A}_0 z_{0} (t)\quad t>0,\\
z_{0}(0)=z.
       \end{array}
       \right.
\end{equation}

Many problems in mathematical physics are led to the study of homogenization of differential operators. There are many works in this direction, we refer to the book of Cioranescu and Donato \cite{CD00}, Zhikov, Kozlov and Oleinik \cite{JKO94} and the book of C. Evans \cite{Evans88} for a list of literature. 

One of the important methods in dealing with periodic homogenization problems is the well-known two scale convergence method, which  is a speciﬁc type of weak convergence developed especially for the theory of homogenization (see, e.g. \cite{Allaire92,Allaire93,Nechvatal03,Ng89,Ng90,ZP04}) in early 90-s in order to overcome the problems connected with weak convergence. This leads to several interesting phenomena in the homogenization process.

In this paper, we are concerned with a  weakly thermoelastic wave model with periodic coefficients. We mean by weakly thermoelastic system, system \eqref{thermoe0} where the coupling terms ($\gamma\theta_x$ and $\gamma u_{tx}$) are replaced by ($\gamma \theta$ and $\gamma u_t$). With these coupling terms, the thermoelastic system considered is different to the one intensively treated in literature, in the sense that its solutions do not decay exponentially to the equilibrium, but pollynomially. A fact that has been shown by
Khodja $\&$ al. in \cite{KBT1997} and has been studied by Liu and Rao in \cite{LR05}. Our purpose in this work is to investigate the problem of homogenization of solutions associated to the weakly coupled thermoelastic wave system by using two scale convergence method \cite{ZP04} and the Trotter Kato convergence theorem \cite{P05,ZP07}. Afterwards, we illustrate numerically the mathematical results.

We end the introduction with some words on the structure of this article. In Section 2, we start by some notations and classical spaces that we will use throughout this work. In Section 3, we recall some properties and results that will be useful for the study of the homogenisation problem. In Section 4, we prove a well posedness result. In Section 5,  we prove some a priori estimates and define some functional spaces for two-scale limits. In Section 6 and 7 we prove the main theorem of our paper. In Section 8, we show that the solutions of our system decays polynomially and not exponentially to zero. Finally, in Section 9, we explain by numerical experiments the difference between exponential and polynomial stability for linear thermoelastic wave systems, and this, in terms of spectral properties
and energy behavior and we show that these results are in perfect agreement with the theory already established in the works \cite{LZb1994,MN16}.

\section{\bf Notations}
Throughout this study, let $\Omega$ be a bounded domain in $\mathbb{R}^N$ with smooth boundary $\partial\Omega$ and $Y = (0,1)^N$ be the periodicity cell with Lebesgue measure $dy$. The symbol $\#$ is used for denotation $Y$-periodic functions. A function $v$ is said to be $Y$-periodic, if $v(y+k)= v(y)$, $y \in Y$, $k\in\Z^N$. If the function $v$ has more variables we say that is $Y$-periodic in $y$. 

We use standard notations. Given a closed linear operator $A$ we denote by $\sigma(A)$, $\rho(A)$, $D(A)$, $rg(A)$ and $R(\lambda,A) = (\lambda-A)^{-1}$ the spectrum of $A$, the resolvent set of $A$, the domain of $A$, the range of $A$ and the resolvent of $A$ respectively. 
\begin{itemize}
\item $Q_T=\Omega\times (0,T)$, $T>0$.
\item $\partial Q_T=\partial\Omega\times (0,T)$, $T>0$.
\item $\C_{-}=\{\lambda\in \C\;:\;\Re\lambda<0\}$
\item $C_0^{\infty}(\Omega)$: the linear spaces of all (test) functions which are infinitely differentiable, and respectively compactly supported in domain $\Omega$.

\item $L^2(\Omega)$: the space of measurable functions $f :\Omega\rightarrow\C$ for which $\{\int_\Omega f^2 dx\}^2<\infty$ and with inner product $\langle f,g\rangle=\int_\Omega f\bar{g}dx$ and associated norm $\|f\|_2:=\langle f,f\rangle^{\frac{1}{2}}$.


\item $H^1(\Omega)$: the space consisting of all integrable functions $f :\Omega\rightarrow\R$ whose ﬁrst-order weak derivatives exist and are square integrable, $H^1(\Omega)=\{f\;/\;f\in L^2(\Omega),\nabla f\in (L^2(\Omega))^N\}$.

\item $H^1_0(\Omega)$: the space consisting of functions in $H^1(\Omega)$ which vanich on the boundary of $\Omega$,
$$H^1_0(\Omega)=\{f\in H^1(\Omega);\/\;f_{|\partial\Omega}=0\}.$$
\item $H^{-1}(\Omega)$: the dual space of $H^1_0(\Omega)$.
\item We denote spaces of functions periodic with respect to
the unit cell $Y$ as
\begin{align*}
F_{\#}(Y)=\{u:\;u\in F_{loc}(\R^N),\; u\text{ is }Y-\text{periodic}\},
\end{align*}
where $F$ may be, for example, $C$, $C^{\infty}$, $L^p$, $H^1$ or $H^{-1}$.
\end{itemize}
By $C$, $C_1$ etc. we denote generic constants which may change from line to line.
\section{\bf Preliminaries}

For the convenience of the reader, we recall the definition of two-scale convergence, which was ﬁrst introduced by Nguetseng \cite{Ng89,Ng90} and then expanded by Allaire \cite{Allaire92,Allaire93}. The terms weak and strong two-scale convergence were introduced by V. V. Zhikov.

\begin{Def}
A sequence of functions $\{u_\varepsilon\}$ that is bounded in $L^2(\Omega)$ is said to be weakly two-scale convergent to $u_0(x,y)\in L^2(\Omega\times Y)$ (we write $u_\varepsilon(x)\overset{2}{\rightharpoonup} u_0(x,y)$) if
\begin{align*}
\lim\limits_{\varepsilon\to 0}\int_{\Omega}u_{\varepsilon}(x)\phi(x)b(\frac{x}{\varepsilon})dx=\int_{\Omega}\int_{Y}u_0(x,y)\phi(x)b(y)dydx,
\end{align*}
for any $\phi\in C_0^{\infty}(\Omega)$, $b\in C_{\#}^{\infty}(Y)$.
\end{Def}
In addition to the two-scale convergence, we introduce the notion of what is sometimes called strong two-scale convergence. Thisconcept is needed to pass to the limit for some products of two-scale convergent sequences.
\begin{Def}
A sequence of functions $\{u_\varepsilon\}$ that is bounded in $L^2(\Omega,dx)$ is strongly two-scale convergent to $u_0(x,y)\in L^2(\Omega\times Y)$ (we write $u_\varepsilon(x)\overset{2}{\to} u_0(x,y)$) if
\begin{align*}
\lim\limits_{\varepsilon\to 0}\int_{\Omega}v_\varepsilon(x) u_{\varepsilon}(x)dx=\int_{\Omega}\int_{Y}v(x,y)u_0(x,y)dydx,\quad\text{for }v_\varepsilon(x)\overset{2}{\rightharpoonup}v(x,y).
\end{align*}
\end{Def}
The following theorem is a modification of Trotter-Kato theorem for variable Banach spaces cf. \cite[Thm. 4.4,p.87]{Pazy83} or \cite[Thm. 4.9,p.212]{EN2000}, to our setting.
\begin{thm}[Weak two-scale Trotter Kato theorem]
Let $\Tbf_\varepsilon (\cdot)$, $\varepsilon>0$, be a contraction $C_0$-semigroup on the Hilbert space $H$, with generator $\mathbf{A}_\varepsilon$, $\varepsilon>0$. Suppose that for any $\lambda > 0$, we have
\begin{itemize}
\item[(i)] weak two-scale resolvent convergence: $R(\lambda,\Abf_\varepsilon)f \overset{2}{\rightharpoonup} R(\lambda)f$
for some $f \in H$ and
\item[(ii)] the range of $R(\lambda)$ is dense in $H$.
\end{itemize}
Then there exists a unique operator $\Abf_0$ generating a contraction $C_0$-semigroup $\Tbf_0(\cdot)$ on the Hilbert space $H$, such that $R(\lambda):= R(\lambda,\Abf_0)$. Moreover,
\begin{align*}
\Tbf_\varepsilon (t)f \overset{2}{\rightharpoonup} \Tbf_0(t)f
\end{align*}
for all $t>0$ and $f\in H$.
\end{thm}
\section{\bf Well-posedness}
In the following we consider the weakly coupled thermoelastic wave model in a bounded domain $\Omega$ with rapidly oscillating coefficients:
\begin{equation}
\label{thermoe}
 \left\{\begin{array}{llll}
                      c(\frac{x}{\varepsilon})\dfrac{\partial^2 u_\varepsilon}{\partial t^2}=div(a(\frac{x}{\varepsilon})\nabla u_\varepsilon)-\gamma(\frac{x}{\varepsilon})\theta_\varepsilon\quad &in\;Q_T\\
                      d(\frac{x}{\varepsilon})\dfrac{\partial\theta_\varepsilon}{\partial t}=div(b(\frac{x}{\varepsilon})\nabla \theta_\varepsilon)
		      +\gamma(\frac{x}{\varepsilon}) \dfrac{\partial u_\varepsilon}{\partial t}\quad &in\;Q_T\\
u_\varepsilon=0=\theta_\varepsilon\quad &on\;\partial Q_T\\ 
u_\varepsilon(x,0)=\tilde{u}_{0}(x),\; \dfrac{\partial u_\varepsilon}{\partial t}(x,0)=\tilde{u}_{1}(x),\; \theta_\varepsilon(x,0)=\tilde{\theta}_{0}(x) \quad &in \;\Omega.
                     \end{array}
                   \right.
\end{equation}
Here $u_\varepsilon$ represents the displacement (longitudinal or transverse, depending upon the application) field and $\theta_\varepsilon$ the temperature increment field. The first equation is the scalar version of the equations of motion and the second is the heat conduction equation. The thermomechanical coupling $\gamma(\cdot)$ results from
consideration of the interaction between deformation and temperature: a temperature change induces strain and conversely. 
The matrices $a(y)$ and $b(y)$, $y\in Y$ are respectively the diffusivity and propagation speed in the medium. 
For some $\alpha\geqslant 1$, we assume that the matrix $a$ satisfyes
\begin{align}
\label{hypo-ineq}
\alpha|\xi|^2\leq \xi^T a(y)\xi\leq \alpha^{-1}|\xi|^2,\quad\forall \xi\in\R^N.
\end{align}

Define the class of matrix functions:
\begin{align*}
E(\Omega)=E(\alpha,\alpha^{-1},\Omega)=\{a=[a_{ij}(x)]: a\text{ is symmetric and satisfies \eqref{hypo-ineq}}\}.
\end{align*}
The first inequality is nothing but the uniform
ellipticity. Moreover, the functions $c$, $d$ and $\gamma$ are bounded positive $Y$-periodic functions
\begin{itemize}
\item $0<\alpha\le c(y)\le \alpha^{-1}<\infty\qquad\forall y\in Y.$
\item $0<\alpha\le d(y)\le \alpha^{-1}<\infty\qquad\forall y\in Y.$
\item $0<\gamma(y)<1\qquad\forall y\in Y.$
\end{itemize}
We define the energy space $H=H^1_0(\Omega)\times L^{2}(\Omega)\times L^{2}(\Omega)$
with the inner product
\[
\left\langle (u,v,\theta),(\tilde{u},\tilde{v},\tilde{\theta})\right\rangle_\varepsilon=\langle A_\varepsilon u,\tilde{u}\rangle_2+\langle c_\varepsilon v,\tilde{v}\rangle_2+\langle d_\varepsilon \theta,\tilde{\theta}\rangle_2
\]
and the induced norm
$$
\|(u,v,\theta)\|_\varepsilon^{2}=\int_\Omega|\big(a_\varepsilon(x)\big)^\frac{1}{2}\nabla u(x)|^{2}+c_\epsilon (x)|v(x)|^{2}+d_\varepsilon (x)|\theta(x)|^{2}dx.
$$
\begin{rmk}
Notice that in the energy norm $\big(a_\varepsilon(y)\big)^\frac{1}{2}$ is well-defined as a square root of a symmetric nonnegative square matrix $a(y)$.
\end{rmk}
The energy associated to the evolution equation \eqref{thermoe} is given by
$$
E_\varepsilon(t)=\frac{1}{2} \|U\|^2_\varepsilon.
$$
\begin{rmk}
\label{rmk4.2}
In view of the properties of the coefficients, $\|\cdot\|_\varepsilon$ is a norm on $H$, equivalent to the natural Sobolev norm on $H$, noted $\|\cdot\|$, that is, if $U$ is in $H$,
\[
\alpha\|U\|^2\leq \|U\|^2_\varepsilon\le\alpha^{-1}\|U\|^2.
\]
\end{rmk}
By introducing new variable (velocity) $v_\varepsilon=\frac{\partial u_\varepsilon}{\partial t}$, the evolution equation can be reduced to the following Cauchy problem
\begin{equation}
\label{cp}
 \left\{\begin{array}{ll}
 \frac{dU_\varepsilon}{dt}=\mathbf{A}_\varepsilon U_\varepsilon,\\
 U_\varepsilon(0)=\tilde{U}_{0}=(\tilde{u}_{0},\tilde{u}_{1},\tilde{\theta}_{0})^{T},\\
\end{array}
\right.
\end{equation}
with
\begin{equation}
 U_\varepsilon=\left(\begin{smallmatrix}
	u_\varepsilon\\
	v_\varepsilon\\
	\theta_\varepsilon
      \end{smallmatrix}
\right)
\end{equation}
and 
\begin{equation}
\mathbf{A}_\varepsilon=\left(\begin{matrix}
                    0 & I & 0\\
	       -\frac{1}{c_\varepsilon} A_{\varepsilon} & 0 & -\frac{\gamma_\varepsilon}{c_\varepsilon} I\\
		    0 & \frac{\gamma_\varepsilon}{d_\varepsilon} I & -\frac{1}{d_\varepsilon} B_{\varepsilon}
                   \end{matrix}
\right).
\end{equation}
Here we have used the notation $A_{\varepsilon}$ and $B_{\varepsilon}$ defined by 
\begin{align*}
A_{\varepsilon}w=-div(a(\frac{\cdot}{\varepsilon})\nabla w),&\qquad B_{\varepsilon}w=-div(b(\frac{\cdot}{\varepsilon})\nabla w).
\end{align*}
The domain of the operator $\Abf_\varepsilon$ is defined by
\[
D(\Abf_\varepsilon)=D(A_\varepsilon)\times D(A_\varepsilon^{\frac{1}{2}})\times D(B_\varepsilon)
\]
with
\[
\Abf_\varepsilon U=\Big(v,-\frac{1}{c_\varepsilon}A_\varepsilon u-\frac{\gamma_\varepsilon}{c_\varepsilon}\theta,\frac{\gamma_\varepsilon}{d_\varepsilon} v-\frac{1}{d_\varepsilon}B_\varepsilon\theta\Big)^T.
\]
\begin{thm}
\label{thm1}
The family of operators $\Abf_\varepsilon$ generates a family of contraction semigroups $\Tbf_\varepsilon (\cdot)$ on the Hilbert space $H$ endowed with the norm $\|\cdot\|_\varepsilon$.
\end{thm}
\begin{dem}
We first prove that $\Abf_\varepsilon$ is a dissipative operator. Indeed, for any $U\in D(\Abf_\varepsilon)$, we have
\begin{align}
\left\langle\Abf_\varepsilon U,U\right\rangle_\varepsilon &=\langle A_\varepsilon v,u\rangle_2-\langle A_\varepsilon u+\gamma_\varepsilon\theta,v\rangle_2+\langle\gamma_\varepsilon v-B_\varepsilon\theta,\theta\rangle_2\nonumber\\
					&=\langle A_\varepsilon v,u\rangle_2-\langle A_\varepsilon u, v\rangle_2-\langle\gamma_\varepsilon\theta,v\rangle_2+\langle\gamma_\varepsilon v,\theta\rangle_2-\langle B_\varepsilon\theta,\theta\rangle_2\nonumber\\
\Re\left\langle\Abf_\varepsilon U,U\right\rangle_\varepsilon &=-\langle B_\varepsilon\theta,\theta\rangle_2\nonumber\\
&\leqslant -\alpha\|\nabla\theta\|^2_2\nonumber.
\end{align}
This implies that $\Abf_\varepsilon$ is a dissipative operator (due to the assumptions on the coefficients).\\
Let's show now that $0\in\rho(\Abf_\varepsilon)$. That is, for all $(f,g,h)^{T}\in H$, there exist $U=(u,v,\theta)^T\in D(\Abf_\varepsilon)$ such that
\begin{align*}
\Abf_\varepsilon U=F,
\end{align*}
i.e.,
\begin{align*}
v &=f\\
\frac{1}{c_\varepsilon}A_\varepsilon u-\gamma_\varepsilon\theta &=g\\
\gamma_\varepsilon v+\frac{1}{d_\varepsilon}B_\varepsilon\theta &=h.
\end{align*}
Taking $v = f$ in the last equation, we obtain $\frac{1}{d_\varepsilon}B_\varepsilon\theta=h-\gamma_\varepsilon f$. Since $h-\gamma_\varepsilon f\in L^2(\Omega)$ and $0\in\rho(B_\varepsilon)$ there exists only one function $\theta\in D(B_\varepsilon)$ that satisfies $\frac{1}{d_\varepsilon}B_\varepsilon\theta=h-\gamma_\varepsilon f$. Similarly 
\[
\frac{1}{c_\varepsilon}A_\varepsilon u=g+\gamma_\varepsilon\theta.
\]
Since $g+\gamma_\varepsilon\theta\in L^2(\Omega)$ and $0\in\rho(A_\varepsilon)$, then there exists only one function $u\in D(A_\varepsilon)$ that satisfies $\frac{1}{c_\varepsilon}A_\varepsilon u=g+\gamma_\varepsilon\theta$. Therefore $U=(u,v,\theta)\in D(\Abf_\varepsilon)$, for all $\varepsilon>0$. Moreover,
\[
\|\Tbf_\varepsilon(t)U\|_\varepsilon\leqslant\|U\|_\varepsilon\quad\text{ for any $U$ in $H$}.
\]
The conclusion immediately follows by the Lumer-Phillips theorem, see \cite[Theorem 4.6., p.16]{Pazy83}.
\end{dem}
\section{\bf A priori estimates and functional spaces for two-scale limits}
In this subsection In this subsection, we follow the procedure of I.V. Kamtoski and V.P. Smyshlyaev \cite{KS2018} to determine the two-scale homogenized limit. For a fixed $\lambda > 0$, we derive in a standard way a priori estimates for the solution $U_\varepsilon$ of the Cauchy problem \eqref{cp}.
\subsection{A priori estimates}
\begin{lem}
\label{lem1}
For $\varepsilon>0$, the following a priori estimates hold:
\begin{align}
\|c_\varepsilon^\frac{1}{2} u_\varepsilon\|_2 &\leq C\|F\|_2\\
\|a_\varepsilon^\frac{1}{2} \nabla u_\varepsilon\|_2 &\leq C\|F\|_2\\
\|d_\varepsilon^\frac{1}{2} \theta_\varepsilon\|_2 &\leq C\|F\|_2\\
\|b_\varepsilon^\frac{1}{2} \nabla \theta_\varepsilon\|_2 &\leq C\|F\|_2
\end{align}
with a constant $C$ independent of $\varepsilon$.
\end{lem}

\begin{dem}
According to Theorem \ref{thm1} and \cite[Chap 9, p.241]{Yosida80} the right half complex plane belongs to the resolvent set of $\Abf_\varepsilon$, for every $\varepsilon$. Let us consider $F=(f,g,h)^T$ to be an element of $H$. We take $\lambda>0$ to be real strictly positive. The following string of equivalences holds: $R(\lambda,\Abf_\varepsilon)F=U_\varepsilon,\quad(U_\varepsilon=(u_\varepsilon,v_\varepsilon,\theta_\varepsilon)^T \in D(\Abf_\varepsilon))$ is equivalent to
\begin{align*}
\lambda u_\varepsilon-v_\varepsilon &=f\\
-\frac{1}{c_\varepsilon}A_\varepsilon u_\varepsilon+\lambda v_\varepsilon+\frac{\gamma_\varepsilon}{c_\varepsilon}\theta_\varepsilon &=g\\
-\frac{\gamma_\varepsilon}{d_\varepsilon} v_\varepsilon+\lambda\theta_\varepsilon-\frac{1}{d_\varepsilon}B_\varepsilon\theta_\varepsilon &=h
\end{align*}
i.e.,
\begin{align*}
v_\varepsilon &=\lambda u_\varepsilon-f\\
(\lambda^2 c_\varepsilon-A_\varepsilon) u_\varepsilon+\gamma_\varepsilon\theta_\varepsilon &=c_\varepsilon g+\lambda c_\varepsilon f\\
-\lambda\gamma_\varepsilon u_\varepsilon+(\lambda d_\varepsilon-B_\varepsilon)\theta_\varepsilon &=d_\varepsilon h-\gamma_\varepsilon f.
\end{align*}
Multiplying the last two equations respectively by $\lambda u_\varepsilon$ and $\theta_\varepsilon$ results in
\begin{align*}
\lambda^3\langle c_\varepsilon u_\varepsilon,u_\varepsilon\rangle-\lambda\langle A_\varepsilon u_\varepsilon, u_\varepsilon\rangle+\lambda\langle \gamma_\varepsilon\theta_\varepsilon,u_\varepsilon\rangle & =\langle c_\varepsilon g+\lambda c_\varepsilon f,\lambda u_\varepsilon\rangle\\
\langle -\lambda\gamma_\varepsilon u\varepsilon,\theta_\varepsilon\rangle+\lambda\langle d_\varepsilon\theta_\varepsilon,\theta_\varepsilon\rangle-\langle B_\varepsilon\theta_\varepsilon,\theta_\varepsilon\rangle & =\langle d_\varepsilon h-\gamma_\varepsilon f,\theta_\varepsilon\rangle
\end{align*}
then summing
\begin{align*}
\lambda^3\|c_\varepsilon^\frac{1}{2} u_\varepsilon\|^2_2+\lambda\|a_\varepsilon^\frac{1}{2}\nabla u_\varepsilon\|^2_2+\lambda\|d_\varepsilon^\frac{1}{2} \theta_\varepsilon\|^2_2+\|b_\varepsilon^\frac{1}{2} \nabla\theta_\varepsilon\|^2_2=\langle c_\varepsilon g+\lambda c_\varepsilon f,\lambda u_\varepsilon\rangle_2+\langle d_\varepsilon h-\gamma_\varepsilon f,\theta_\varepsilon\rangle_2.
\end{align*}
In view of the properties of the coefficients and using Young's inequality and Poincaré inequality
\begin{align*}
\lambda^3\|c_\varepsilon^\frac{1}{2} u_\varepsilon\|^2_2+\lambda\|a_\varepsilon^\frac{1}{2}\nabla u_\varepsilon\|^2_2+\lambda\|d_\varepsilon^\frac{1}{2} \theta_\varepsilon\|^2_2+\|b_\varepsilon^\frac{1}{2} \nabla\theta_\varepsilon\|^2_2\leq C\|F\|_2^2.
\end{align*}
\end{dem}

\subsection{Functional spaces for two-scale limits}

We next also introduce the following 'dual' spaces $W_a$ and $W_b$ of admissible 'microscopic fluxes', of tensor fields on $Y$:
\begin{align}
\tag{**}
W_a:&=\Bigg\{\psi\in \Big(L^2_{\#}(Y)\Big)^N \Big| div_y\Bigg(\Big((a(y)\Big)^\frac{1}{2}\psi(y)\Bigg)=0\quad\text{in } H^{-1}_{\#}(Y) \Bigg\}\\
\tag{***}
W_b:&=\Bigg\{\psi\in \Big(L^2_{\#}(Y)\Big)^N \Big| div_y\Bigg(\Big((b(y)\Big)^\frac{1}{2}\psi(y)\Bigg)=0\quad\text{in } H^{-1}_{\#}(Y) \Bigg\}
\end{align}

It immediately follows from the definitions above that $W_a$ and $W_b$ are closed linear subspaces of the Hilbert space $\Big(L^2_{\#}(Y)\Big)^N$, and hence can themselves be regarded as Hilbert spaces with respective inherited $L^2$ inner product. When there is no confusion we denote $W_a$ and $W_b$ by $W$.

We will additionally introduce, in a standard way, Hilbert spaces $L^2(\Omega; H^1_{\#}(Y))$, $L^2(\Omega;W_a)$ and $L^2(\Omega;W_b)$ of functions of two independent variables $x\in\Omega$ and $y\in Y$, which can thereby be regarded as functions of $x$ with values in the appropriate (Hilbert) space.

The a priori estimates in Lemma \ref{lem1}, via adapting accordingly the properties of the two-scale convergence, imply the following
\begin{lem}
\label{lem2}
There exist $u_0(x,y),\theta_0(x,y)\in L^2(\Omega)$ and $\xi_0 (x,y)\in L^2(\Omega;W_a)$ and $\chi_0(x,y)\in L^2(\Omega;W_b)$ such that, up to extracting a subsequence in $\varepsilon$ which we do not relabel,
\begin{align}
u_\varepsilon & \overset{2}{\rightharpoonup} u_0(x,y)\\
\varepsilon \nabla u_\varepsilon & \overset{2}{\rightharpoonup} \nabla_y u_0(x,y)\\
\Big((a(x/\varepsilon)\Big)^\frac{1}{2}\nabla u_\varepsilon & \overset{2}{\rightharpoonup} \xi_0(x,y)\\
\theta_\varepsilon & \overset{2}{\rightharpoonup} \theta_0(x,y)\\
\varepsilon \nabla \theta_\varepsilon & \overset{2}{\rightharpoonup} \nabla_y \theta_0(x,y)\\
\Big((b(x/\varepsilon)\Big)^\frac{1}{2}\nabla \theta_\varepsilon & \overset{2}{\rightharpoonup} \chi_0(x,y)
\end{align}
\end{lem}

\begin{dem}
\begin{enumerate}
\item[(1)] According to the theorem on (weak) two-scale compactness of a bounded sequence in $L^2(\Omega)$, the á priori estimate (5.2) implies, up to extracting a subsequence in $\varepsilon$ (not relabelled), the
existence of a weak two-scale limit $\xi_0\in\Bigg(L^2\Big(\Omega\times Y\Big)\Bigg)^N=L^2\Bigg(\Omega;\Big(L^2_{\#}(Y)\Big)^N\Bigg)$, which yields (5.7).

We show that in fact $\xi_0(x, y)\in L^2(\Omega;W_a)$. By definition of (weak) two scale convergence, for $\varphi\in C_0^\infty(\Omega)$ and $b\in C_{\#}^\infty(Y)$
\begin{align*}
\lim_{\varepsilon\to 0}&\int_\Omega a(x/\varepsilon)\nabla u_\varepsilon.\varepsilon\nabla\Big(\varphi(x)b(x/\varepsilon)\Big)dx\\
&=\int_\Omega\varphi(x)\int_Y\Big(a(y)\Big)^\frac{1}{2}\xi_0(x,y).\nabla_y b(y)dy dx=0,
\end{align*}
where we have used the assumption on the coefficients $a_\varepsilon$. The density of $\varphi(x)$ in $L^2(\Omega)$ implies that for all $b\in C_{\#}^\infty(Y)$ the inner integral is zero for a.e. $x\in\Omega$. Since $b(y)$ are in turn dense in $H^1_{\#}(Y)$, this implies that, for a.e. $x$, $\xi_0(x,\cdot)$ obeys (**) and hence $\xi_0(x,\cdot)\in W_a$ implying $\xi_0\in L^2(\Omega;W_a)$
\item[(2)] Further, according to e.g. \cite[Prop. 1.14 (ii)]{Allaire92}, (5.1) together with (5.2) imply (5.5)–(5.6) for some $u_0(x,y)\in L^2(\Omega; H^1_{\#}(Y))$.

Show finally that in fact $u_0(x, y)\in L^2(\Omega)$. For any $\varphi\in C_0^\infty(\Omega)$ and $b\in C_{\#}^\infty(Y)^N$,
\begin{align}
\lim_{\varepsilon\to 0}&\int_\Omega \Big(a(x/\varepsilon)\Big)^\frac{1}{2}\varepsilon\nabla u_\varepsilon(x).\varphi(x)b(x/\varepsilon)dx=\int_\Omega\int_Y\Big(a(y)\Big)^\frac{1}{2}\nabla_y u_0(x,y)\varphi(x)b(y)dy dx,
\end{align}
where we have used (5.6).

On the other hand, (5.2) ensures that
\begin{align*}
\Bigg\|\Big(a(x/\varepsilon)\Big)^\frac{1}{2}\varepsilon\nabla u_\varepsilon (x)\Bigg\|_2\to 0,
\end{align*}
and hence the limit in (5.11) is zero. This implies for the right hand side of (5.11),
\begin{align*}
\int_\Omega\varphi(x)\int_Y\Big(a(y)\Big)^\frac{1}{2}\nabla_y u_0(x,y).b(y)dy dx=0,\quad\forall\varphi\in C_0^\infty(\Omega),\: b\in C_{\#}^\infty(Y)^N.
\end{align*}
By density of $\varphi$ and $b$, this gives
\begin{align}
\Big(a(y)\Big)^\frac{1}{2}\nabla_y u_0(x,y)=0\text{ for a.e. }x,
\end{align}
and therefore, pre-multiplying (5.12) by $\Big(a(y)\Big)^\frac{1}{2}$, yields $u_0(x,y)\in L^2(\Omega)$.
\item[(3)] The remaining (weak) two scale convergence (5.8), (5.9) and (5.10) can be shown in a similar way.
\end{enumerate}
\end{dem}

\begin{rmk}
According to estimates (5.6) with (5.2) we can show that $u_0(x,y)=u_0(x)$. Similarly (5.4) and (5.9) yield $\theta_0(x,y)=\theta_0(x)$.
\end{rmk}

\section{\bf The two-scale limit problem}
We establish first an important property connecting, the generalized
two-scale limit flux $\xi_0(x, y)$ to the two-scale limit field $u_0(x, y)$, see Lemma \ref{lem2}.

We introduce the following set of 'product' test functions in $L^2(\Omega;W)$. Let $\Psi(x, y) = \varphi(x)b(y)$, where $\varphi\in C_c^{\infty}(\overline{\Omega})$ and $b\in W$. Here $C_c^{\infty}(\overline{\Omega})$ consists of restrictions to $\Omega$ of all the scalar functions in $\Omega$ from $C_0^{\infty}(\R^N)$, i.e. of infinitely differentiable functions with a compact support
in the whole of $\R^N$. We note that the linear span of such test functions $\Psi$ is dense in $L^2(\Omega;W)$: e.g. an arbitrary $\Psi(x, y)\in L^2(\Omega;W)\subset\Big
(L^2(\Omega\times Y)\Big)^N$ is approximated in $\Big
(L^2(\Omega\times Y)\Big)^N$ by
linear span of $\varphi_m(x)\tilde{b}_m(y)$ with $\varphi_m(x)\in C_0^{\infty}(\Omega)\subset C_c^{\infty}(\overline{\Omega})$ and $\tilde{b}_m(y)\in C_{\#}^{\infty}(Y)$, and then setting $b_m(y) = P_W \tilde{b}_m(y)$ with $P_W$ denoting orthogonal projection on $W$ in $\Big(L^2_{\#}(Y)\Big)^N$.

The following important lemma holds.

\begin{lem}
Let $u_0(x,y)$ and $\xi_0(x,y)$ be as in Lemma \ref{lem2}. Then the following integral identity holds
\begin{align}
\forall\varphi\in C_c^\infty(\overline{\Omega}),\: b\in W_a,\quad &\int_\Omega\int_Y \xi_0(x,y)\varphi(x)b(y)dydx=\nonumber\\
& \int_\Omega\int_Y u_0(x,y).div_x\Bigg((a(y))^\frac{1}{2}\varphi(x)b(y)\Bigg)dydx.
\end{align}
\end{lem}

\begin{dem}
Let $\varphi\in C_0^\infty(\Omega)$ and $b\in W_a$. Then, by (5.7)
\begin{align}
\lim_{\varepsilon\to 0}\int_\Omega \Big(a(x/\varepsilon)\Big)^\frac{1}{2}\nabla u_\varepsilon(x).\varphi(x)b(x/\varepsilon)dx=\int_\Omega\int_Y\xi_0(x,y)\varphi(x)b(y)dy dx.
\end{align}
On the other hand, integrating by parts and using $u_\varepsilon\in H^1_0(\Omega)$ and (**), gives
\begin{align*}
\int_\Omega \Big(a(x/\varepsilon)\Big)^\frac{1}{2}\nabla u_\varepsilon(x).\varphi(x)b(x/\varepsilon)dx=-\int_\Omega u_\varepsilon(x).div_x\Bigg(\Big(a(y)\Big)^\frac{1}{2}\varphi(x)b(y)\Bigg)\Bigg|_{y=x/\varepsilon}dx.
\end{align*}
Passing to the limit and using (5.5) then yields
\begin{align}
\lim_{\varepsilon\to 0}\int_\Omega \Big(a(x/\varepsilon)\Big)^\frac{1}{2}\nabla u_\varepsilon(x).\varphi(x)b(x/\varepsilon)dx=\int_\Omega\int_Y u_0(x,y).div_x\Bigg(\Big(a(y)\Big)^\frac{1}{2}\varphi(x)b(y)\Bigg)dxdy.
\end{align}
Comparing (5.14) and (5.15) results in identity (5.13).
\end{dem}

\begin{rmk}
Notice that the test function $\varphi(x)$ in the above lemma does not have to vanish when $x\in\partial\Omega$. Formally, if we use integration by parts on the right hand side of (5.13), with $b\in W$ we find
\begin{align*}
\int_\Omega\int_Y \Big[\xi_0(x,y)-\Big(a(y)\Big)^\frac{1}{2}\nabla_x u_0(x,y)\Big]\varphi(x)b(y)dydx=0
\end{align*}
i.e., $\xi_0(x,y)-\Big(a(y)\Big)^\frac{1}{2}\nabla_x u_0(x,y)\perp\varphi(x)b(y)$ in $L^2(\Omega\times Y)^N$. Recall that the orthogonal of divergence-free functions are exactly the gradients (see, if necessary, \cite{Tartar78} or \cite{Temam79}). Thus, we deduce that
there exists a unique function $u_1(x,y)$ in $L^2(\Omega;H^1_{\#}(Y))$ solving the cell problem $-div_y\Big(a(y)[\nabla_y u_1(x,y)+\nabla_x u_0(x)]\Big)=0$, such that 
\begin{align*}
\xi_0(x,y)=(a(y))^\frac{1}{2}[\nabla_y u_1(x,y)+\nabla_x u_0(x)].
\end{align*}
Similarly we can show that there exists a unique function $\theta_1(x,y)$ in $L^2(\Omega;H^1_{\#}(Y))$ solving the cell problem $-div_y\Big(b(y)[\nabla_y \theta_1(x,y)+\nabla_x \theta_0(x)]\Big)=0$, such that 
\begin{align*}
\chi_0(x,y)=(b(y))^\frac{1}{2}[\nabla_y \theta_1(x,y)+\nabla_x \theta_0(x)].
\end{align*}
\end{rmk}

\section{\bf The weak two-scale resolvent convergence}

\begin{thm}
\begin{itemize}
\item[(i)] The strongly continuous contraction semigroups $\Tbf_\varepsilon(\cdot)$ associated with $\Abf_\varepsilon$ weakly two-scale converge to the semigroup $\Tbf_0(\cdot)$ associated with $\Abf_0$, i.e. if $f\in H$ then for all $t\geq 0$,
\[
\Tbf_\varepsilon(t)f\overset{2}{\rightharpoonup}\Tbf_0(t)f.
\]
\item[(ii)] For $\varepsilon>0$, the family of Cauchy problems
\begin{align}
\frac{dU_\varepsilon}{dt}=\mathbf{A}_\varepsilon U_\varepsilon,\;U_\varepsilon(0)=\tilde{U}_{0}=(\tilde{u}_{0},\tilde{u}_{1},\tilde{\theta}_{0})^{T},
\end{align}
(weakly) two-scale converge to the limit Cauchy problem
\begin{align}
\frac{dU_0}{dt}=\Abf_0 U_0,\;U_0(0)=\tilde{U}_{0}=(\tilde{u}_{0},\tilde{u}_{1},\tilde{\theta}_{0})^{T},
\end{align}
\end{itemize}
\end{thm}

\begin{dem}
\begin{itemize}
\item[(i)] Let us consider $F=(f,g,h)^T$ to be an element of $H$. We take $\lambda>0$ to be real strictly positive. The following string of equivalences holds: $R(\lambda,\Abf_\varepsilon)F=U_\varepsilon,\quad(U_\varepsilon=(u_\varepsilon,v_\varepsilon,\theta_\varepsilon)^T \in D(\Abf_\varepsilon))$ is equivalent to
\begin{align*}
\lambda u_\varepsilon-v_\varepsilon &=f\\
-\frac{1}{c_\varepsilon}A_\varepsilon u_\varepsilon+\lambda v_\varepsilon+\frac{\gamma_\varepsilon}{c_\varepsilon}\theta_\varepsilon &=g\\
-\frac{\gamma_\varepsilon}{d_\varepsilon} v_\varepsilon+\lambda\theta_\varepsilon-\frac{1}{d_\varepsilon}B_\varepsilon\theta_\varepsilon &=h
\end{align*}
i.e.,
\begin{align*}
v_\varepsilon &=\lambda u_\varepsilon-f\\
(\lambda^2 c_\varepsilon-A_\varepsilon) u_\varepsilon+\gamma_\varepsilon\theta_\varepsilon &=c_\varepsilon g+\lambda c_\varepsilon f\\
-\lambda\gamma_\varepsilon u_\varepsilon+(\lambda d_\varepsilon-B_\varepsilon)\theta_\varepsilon &=d_\varepsilon h-\gamma_\varepsilon f.
\end{align*}
Multiplying respectively by $\phi(x)$, $\varphi(x)$ and $\psi(x)$ in $C_0^{\infty}(\Omega)$ yields
\begin{align*}
&\lambda\int_\Omega u_\varepsilon(x)\phi(x)dx-\int_\Omega v_\varepsilon(x)\phi(x)dx=\int_\Omega f(x)\phi(x)dx,\\
&\lambda^2\int_\Omega c_\varepsilon(x) u_\varepsilon(x)\varphi(x)dx+\int_\Omega a_\varepsilon(x)\nabla u_\varepsilon(x)\nabla\varphi(x)dx+\int_\Omega\gamma_\varepsilon(x)\theta_\varepsilon(x)\varphi(x)dx\\
&=\int_\Omega c_\varepsilon(x)(g(x)+\lambda f(x))\varphi(x)dx\\
&-\lambda\int_\Omega\gamma_\varepsilon(x) u_\varepsilon(x)\psi(x)dx+\lambda\int_\Omega d_\varepsilon(x)\theta_\varepsilon(x)\psi(x)dx+\int_\Omega b_\varepsilon(x)\nabla\theta_\varepsilon(x)\nabla\psi(x)dx\\
&=\int_\Omega d_\varepsilon(x)h(x)\psi(x)dx-\int_\Omega\gamma_\varepsilon(x)f(x)\psi(x)dx.
\end{align*}
Passing to (weak) two-scale limit and using Lemma \ref{lem2}
\begin{align*}
&\lambda\int_\Omega u_0(x)\phi(x)dx-\int_\Omega v_0(x)\phi(x)dx=\int_\Omega f(x)\phi(x)dx,\\
&\lambda^2\int_\Omega<c> u_0(x)\varphi(x)dx+\int_\Omega a^{hom}\nabla u_0(x)\nabla\varphi(x)dx+\int_\Omega\theta_0(x)\varphi(x)dx\\
&=\int_\Omega <c>(g(x)+\lambda f(x))\varphi(x)dx,\\
&-\lambda\int_\Omega<\gamma> u_0(x)\psi(x)dx+\lambda\int_\Omega <d>\theta_0(x)\psi(x)dx+\int_\Omega b^{hom}\nabla\theta_0(x)\nabla\psi(x)dx\\
&=\int_\Omega <d>h(x)\psi(x)dx-\int_\Omega<\gamma>f(x)\psi(x)dx,
\end{align*}
where $<c>=\int_Y c(y)dy$, $<d>=\int_Y c(y)dy$, $<\gamma>=\int_Y\gamma(y)dy$, are respectively the mean value of the functions $c,d$ and $\gamma$. The homogenized matrices $a^{hom}$ and $b^{hom}$ are defined by their entries
\begin{align*}
a^{hom}_{ij} &=\int_Y[(a(y)\nabla_y M_i).e_j+a_{ij}(y)]dy=\int_Y a(y)(e_i+\nabla_y M_i).(e_j+\nabla_y M_j)dy.
\\
b^{hom}_{ij} &=\int_Y[(b(y)\nabla_y N_i).e_j+b_{ij}(y)]dy=\int_Y b(y)(e_i+\nabla_y N_i).(e_j+\nabla_y N_j)dy
\end{align*}
where $N_i$ and $M_i$ are respectively the solutions of the following cell problems:
\begin{align*}
\left\{\begin{array}{ll}
 -div_y a(y)(e_i+\nabla_y M_i(y))=0\quad&\text{in } Y,\\
 y\to M_i(y)\quad &Y\text{-periodic}
\end{array}
\right.
\end{align*}
and
\begin{align*}
\left\{\begin{array}{ll}
 -div_y b(y)(e_i+\nabla_y N_i(y))=0\quad&\text{in } Y,\\
 y\to N_i(y)\quad &Y\text{-periodic}
\end{array}
\right.
\end{align*}
$(e_i)_{1\le i\le N}$ is the canonical basis of $\R^N$. 

Since the range $rg R(\lambda,\Abf_\varepsilon) = D(\Abf_\varepsilon)$ is dense in $H$ and $R(\lambda,\Abf_\varepsilon)f \overset{2}{\rightharpoonup} R(\lambda)f$ for all $f\in H$, we conclude by the weak two-scale Trotter Kato theorem.
\item[(ii)] According to the above resulat, one conclude that the homogenized Cauchy problem is given by  
\[
\frac{dU_0}{dt}=\Abf_0 U_0,\;U_0(0)=\tilde{U}_{0}=(\tilde{u}_{0},\tilde{u}_{1},\tilde{\theta}_{0})^{T},
\]
where the homogenized operator
\begin{equation*}
\mathbf{A}_0=\left(\begin{matrix}
                    0 & I & 0\\
	       -\frac{1}{<c>} A_{hom} & 0 & -\frac{<\gamma>}{<c>} I\\
		    0 & \frac{<\gamma>}{<d>} I & -\frac{1}{<d>} B_{hom}
                   \end{matrix}
\right).
\end{equation*}
Here we have used the notation $A_{hom}$ and $B_{hom}$ defined by 
\begin{align*}
A_{hom}w=-div(a^{hom}\nabla w),&\qquad B_{hom}w=-div(b^{hom}\nabla w)
\end{align*}
$a^{hom}$ and $b^{hom}$ are constant matrices defined as above.

\end{itemize}
\end{dem}

\begin{rmk}
\begin{itemize}
\item[(i)] $M_i(y)$ and $N_i(y)$ are respectively the local variation of displacement and the local variation of temperature created by an averaged (or macroscopic) gradient $e_i$. By linearity, it is not difficult to
compute $u_l (x,y)$, solution of the cell problem $-div_y\Big(a(y)[\nabla_y u_1(x,y)+\nabla_x u_0(x)]\Big)=0$, in terms of $u_0(x)$ and  $M_i(y)$
\begin{align*}
u_1(x,y)=\sum_{i=1}^N\frac{\partial u_0}{\partial x_i}(x)M_i(y).
\end{align*}
Similarly 
\begin{align*}
\theta_1(x,y)=\sum_{i=1}^N\frac{\partial \theta_0}{\partial x_i}(x)N_i(y).
\end{align*}
\item[(ii)] The constant tensors $a^{hom}$ and $b^{hom}$ describe the effective or homogenized properties of the heterogeneous material $a(\frac{x}{\varepsilon})$ and $b(\frac{x}{\varepsilon})$. Note that $a^{hom}$ and $b^{hom}$ do not depend on the choice of domain $\Omega$, functions $c,d$ and $\gamma$, or boundary condition on $\partial\Omega$.
\end{itemize}
\end{rmk}

\section{\bf Asymptotic behaviour of solutions}
We recall here the notions of stability that we encounter in this work.
\begin{Def}
Assume that $A$ is the generator of a $C_0-$semigroup of contractions $(T(t))_{t\geqslant 0}$ on a Hilbert space $H$. The $C_0-$semigroup $T(\cdot)$ is said to be
\begin{enumerate}
\item[1.] strongly stable if
\begin{align*}
\lim\limits_{t \rightarrow +\infty} \|T(t)U_{0}\|_{H}=0,\quad U_{0}\in H;
\end{align*}
\item[2.] exponentially (or uniformly) stable with decay rate $\omega > 0$ if there exists a
constant $M\geqslant 1$ such that 
$$
\|T(t)U_0\|_H\leqslant Me^{-\omega t}\|U_0\|_H,\quad t>0,\; U_0\in H;
$$
\item[3.] polynomially stable (of order $\alpha>0$) if it is bounded, if $i\mathbb{R}\subset\rho(A)$ and if
$$
\|T(t)U_0\|_H\leqslant Ct^{-\alpha}\|A U_0\|_H,\quad t>0,\; U_0\in D(\mathcal{A})
$$
for some positive constant $C$.\\
In that case, one says that solutions of \eqref{cp} decay at a rate $t^{-\alpha}$. The $C_0-$semigroup $T(\cdot)$ is said to be polynomially stable with optimal decay rate $t^{-\alpha}$ (with $\alpha>0$) if it is polynomially stable with decay rate $t^{-\alpha}$ and, for any $\epsilon>0$ small enough, there exists solutions of \eqref{cp} which do not decay at a rate $t^{-(\alpha-\epsilon)}$.
\end{enumerate}
\end{Def}

\begin{thm}
Let us consider $\varepsilon>0$. Then it follows that
\begin{enumerate}
\item[(a)] the sequence of semigroups $\Tbf_\varepsilon(\cdot)$ is strongly stable but not exponentially stable,
\item[(b)] the semigroup  associated with $\Abf_\varepsilon$ is polynomially stable (the order $\alpha=\frac{1}{2}$ is optimal).
\end{enumerate}
\end{thm}

\begin{dem}
The proof follows directly from \cite{LR05,MN16} and Remark \ref{rmk4.2}.
\end{dem}

\section{\bf Numerical Experiments}

\subsection{\bf Approximation of the thermoelastic wave model}
We consider the numerical approximation of the following coupled thermoelastic wave models
\begin{equation}
\label{eq1}
 \left\{\begin{array}{llll}
                      u_{tt}(x,t)- \Delta u(x,t)+\gamma\theta_x(x,t)=0\quad &in\;\Omega\times(0, \infty)\\
		      \theta_{t}(x,t)-\Delta\theta(x,t)-\gamma u_{tx}(x,t) =0\quad &in\;\Omega\times(0, \infty)\\
		      u(x,t)=0=\theta(x,t)\qquad &\text{on } \partial\Omega\times(0,\infty)\\
u(x,0)=u_{0}(x),\; u_{t}(x,0)=u_{1}(x),\; \theta(x,0)=\theta_{0}(x) \quad &on \;\Omega
                     \end{array}
                   \right.
\end{equation}

\begin{equation}
\label{eq2}
 \left\{\begin{array}{llll}
                      u_{tt}(x,t)- \Delta u(x,t)+\gamma\theta(x,t)=0\quad &in\;\Omega\times(0, \infty)\\
		      \theta_{t}(x,t)-\Delta\theta(x,t)-\gamma u_{t}(x,t) =0\quad &in\;\Omega\times(0, \infty)\\
		      u(x,t)=0=\theta(x,t)\qquad &\text{on } \partial\Omega\times(0,\infty)\\
u(x,0)=u_{0}(x),\; u_{t}(x,0)=u_{1}(x),\; \theta(x,0)=\theta_{0}(x) \quad &on \;\Omega
                     \end{array}
                   \right.
\end{equation}
It is well known from literature that system (\ref{eq1}) and (\ref{eq2}) are respectively exponentially and polynomially stable, see \cite{Hansen1992, LZb1993} and \cite{KBT1997, KBB1999, LR05}.

In this section, we will show by numerical experiments, how the coupling terms affect quantitative and qualitative properties of thermoelastic systems $\eqref{eq1}$ and $\eqref{eq2}$. These results could be found in \cite{LZb1994,MN16}.

To do this, we consider a semi discretization version of both systems $\eqref{eq1}$ and $\eqref{eq2}$, obtained with spectral element method. That is, we choose the eigenvectors of the system as the basis vectors. Here, we will use the eigenvectors of the uncoupled thermoelastic system, i.e., $\gamma=0$
in $\eqref{eq1}$ and $\eqref{eq2}$. Let
$$\phi_{j}(x)=\sqrt{\frac{2}{\pi}}\frac{1}{j}\sin jx,\quad \psi_j=\sqrt{\frac{2}{\pi}}\sin jx,\quad \xi_j=\sqrt{\frac{2}{\pi}}\sin jx\;\;j=1,\dots,n.$$
A straightforward calculation yields the following form
\begin{equation}
\label{eq3}
 (S_i)\left\{\begin{array}{ll}
 z'_n(t)=A_{i,n} z_n(t),\quad t\geq 0,\;n\in\mathbb{N},\;i=1,2\\
 z_n(0)=z_{n0},\;n\in\mathbb{N}\\
\end{array}
\right.
\end{equation}
where for all $n\in\mathbb{N}$, $z_n=(u_n,v_n,\theta_n)^T$ is the semi discrete solution, $z_{n0}$ is the discretized initial data, $A_{i,n}$ the discretized dynamic and the subscript $\cdot_i$ refers to system $\eqref{eq1}$ and $\eqref{eq2}$ with  
\begin{align*}
A_{i,n}=\left(\begin{matrix}
                    0 & D_n & 0\\
	       -D_n & 0 & -\gamma F_{i,n}\\
		    0 & \gamma F_{i,n} & -D_n^2
                   \end{matrix}
\right),\quad i=1,2
\end{align*}
and $F_{2,n}=I_n$, 
\begin{align*}
D_n=\left(\begin{matrix}
                    1 &  & \\
	        & \ddots &  \\
		     &   & n
                   \end{matrix}
\right),\quad (F_{1,n})_{ij}=\left\{\begin{array}{ll}
 -\frac{4}{\pi}\frac{ij}{i^2-j^2},\quad |i-j|=\text{odd},\\
 0,\qquad\qquad\text{otherwise}.
\end{array}
\right.
\end{align*}
It has been shown in \cite{LZb1994,MN16} that the approximation semigroups $T_n(t)$ associated to \eqref{eq3} converge strongly to $T(t)$. This has further important implications and we state this in the following theorem.
\begin{thm}
$T_{n}(t),\;T^{*}_{n}(t)\underrightarrow{s} T(t),\;T^{*}(t)$ in $H$, respectively. 
Moreover, the convergence is uniform in bounded $t$-intervals.
\end{thm}
\subsection{\bf Numerical approximation of the spectrum}

Figure 1 and Figure 2 show how the coupling terms affect the placement of eigenvalues of the dynamic $A_{n}$. In Figure 1, we see that a uniform distance between the eigenvalues and the imaginary axis is preserved, see Table 1. Another observation is that for fixed $n$, the eigenvalues of higher frequency modes, in particular, the one of the $n^{th}$ mode, are closer to the imaginary axis. Moreover, as the number of modes increases, these eigenvalues bend back towards the vertical line $\lambda=-\frac{\gamma^2}{2}$, a fact which has been already shown in \cite{KBT1997}. Therefore, the corresponding spectral element approximation scheme preserves the property of exponential stability.
\begin{table}[H]
\begin{center}
\caption{Distance between $\sigma(A_{1,n})$ and the imaginary axis for the spectral element method}
\begin{tabular}{|c|c|} 
\hline
$n$ & min\{-Re\;$\lambda$, $\lambda\in\sigma(A_{1,n})$\}\\
\hline 
8 & 8.9227$\times 10^{-4}$\\
\hline
16 & 8.9383$\times 10^{-4}$\\
\hline
24 & 8.9402$\times 10^{-4}$\\
\hline
32 & 8.9407$\times 10^{-4}$\\ 
\hline
\end{tabular}
\end{center}
\end{table}
\begin{figure}[H]
\label{fig1}
\begin{center}
\includegraphics[scale=0.6]{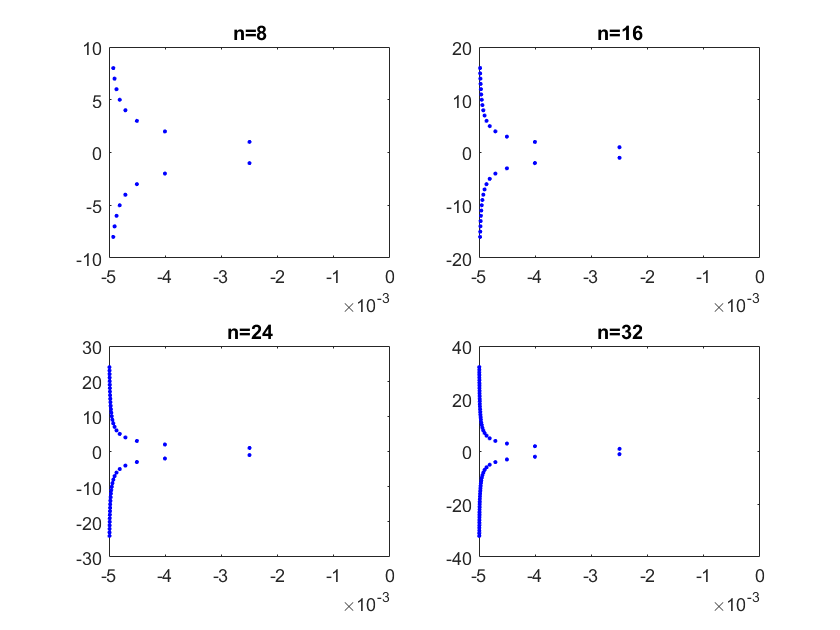}
\caption{Location of the complex eigenvalues of the matrix $A_{1,n}$ with the spectral element method}
\end{center}
\end{figure}
In Figure 2, conversely to Figure 1 where a uniform distance between the eigenvalues and the imaginary axis is preserved, we observe that, as the number of modes increases, an asymptotic behaviour appears in the neighborhood of the imaginary axis at $\pm\infty$. This property is mainly related to systems with polynomial decay, see \cite{BEPS2006}.
\begin{figure}[H]
\begin{center}
\includegraphics[scale=0.5]{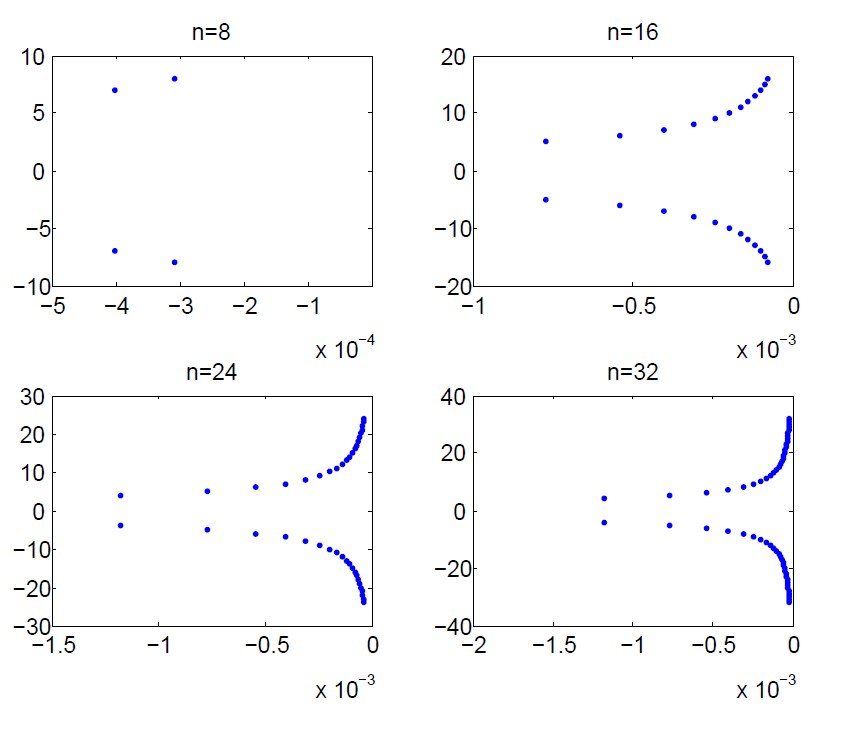}
\caption{Location of the complex eigenvalues of the matrix $A_{2,n}$ with the spectral element method.}
\end{center}
\end{figure}

\subsection{\bf Uniform and polynomial decay of the energy}
The discrete energy associated to system \eqref{eq3} is given by
\begin{equation}
\label{eq4}
E_{i,n}(t)=\frac{1}{2}\sum_{j=1}^n\Big\{|u_j(t)|^{2}+|v_j(t)|^{2}+|\theta_j(t)|^{2}\Big\},\;i=1,2.
\end{equation}
The discrete energy $E_{1,n}$ associated to system $\eqref{eq1}$ decays exponentially to zero, see Figure 3, in the following sense: $\exists M,\alpha$ positive constants such that
\[
E_{1,n}(t)\leqslant Me^{-\alpha t}E_{1,n}(0),\;n\in\mathbb{N},\;t>0.
\]
However, the introduction of the weak coupling term in system $\eqref{eq1}$ has changed the dynamic and consequently the behavior of energy \eqref{eq4}. In this case, we say that system \eqref{eq2} decays polynomially to zero, see Figure 4, in the following sense: $\exists M,\alpha$ positive constants such that
\[
E_{2,n}(t)\leqslant \frac{M}{t}\|A_{2,n}z_{n0}\|^2,\;n\in\mathbb{N},\;t>0.
\]
\begin{figure}[H]
\begin{center}
\includegraphics[scale=0.8]{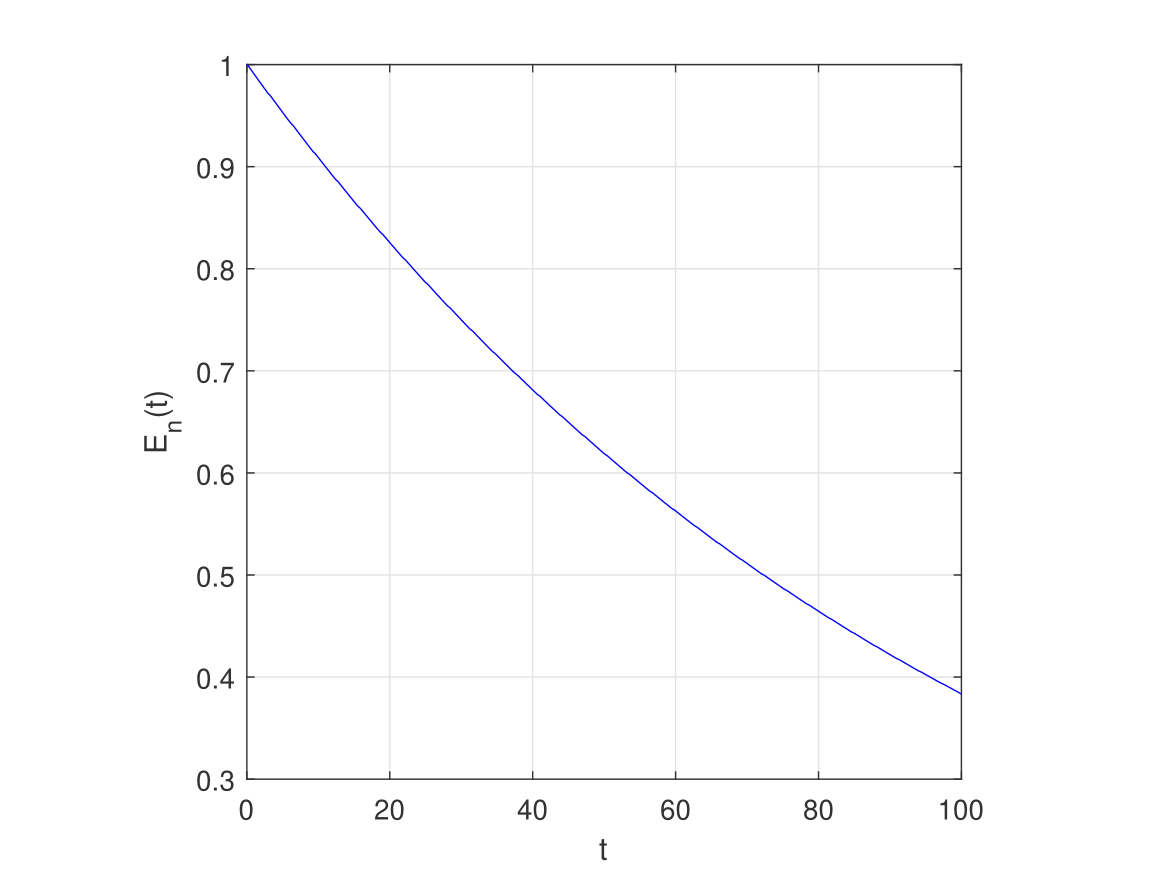}
\caption{Exponential decay of $E_{1,n}(t)$}
\end{center}
\begin{center}
\includegraphics[scale=0.8]{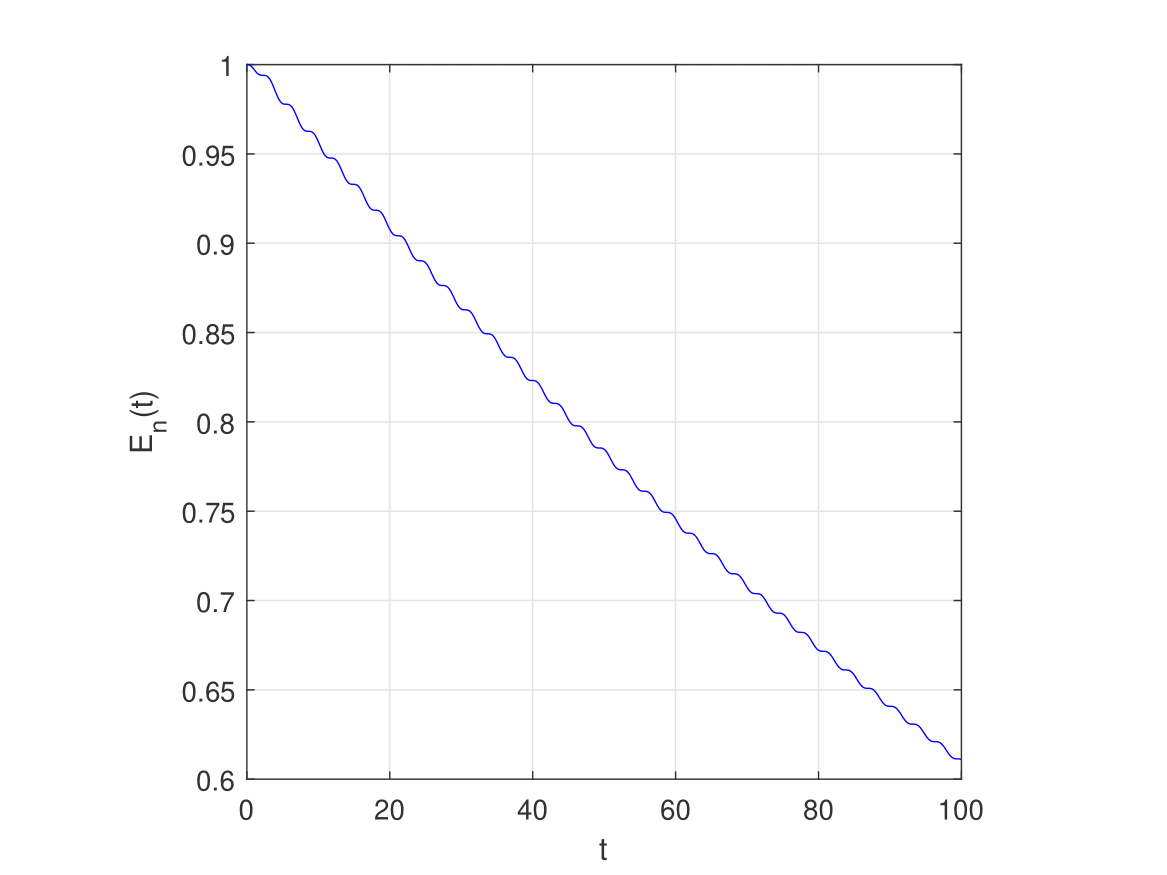}
\caption{Polynomial decay of $E_{2,n}(t)$}
\end{center}
\end{figure}

\subsection{\bf Effect of smoothness of the initial data on the rate of decay of energy} 
It has been shown theoretically, see \cite{BEPS2006, BT2010}, that the energy associated to system \eqref{eq2} is very sensitive to the smoothness of its initial data. This fact, has been also observed numerically, see Figure 6. 
we use a uniform mesh with $n = 100$ elements, fix the final time at $T = 100$, use dt = 0.1 and consider the following initial conditions for $u$, $u_t$ and $\theta$
\[
u(x,0) = 0,\quad \theta(x, 0) = 0,\quad u_t(x, 0) =sin(jx),\; j = 1, 2, 3.
\]
Through Figure 6, we notice that for $j = 1$, the approximate energy $E_{2,n}(t)$ decays to zero as the time $t$ increases. Moreover, we observe that the decay rate depends strongly to j. That is, when $j$ increases, initial data are very oscillating. We say in this case that the rate of decay of the discrete energy $E_{2,n}(t)$ is very sensitive to
the choice of the initial data. However, the behavior of the energy assosiated to system \eqref{eq1} remains indifferent to the smoothness of initial data when $n\to\infty$, see Figure 5.
\begin{figure}[H]
\begin{center}
\includegraphics[scale=0.8]{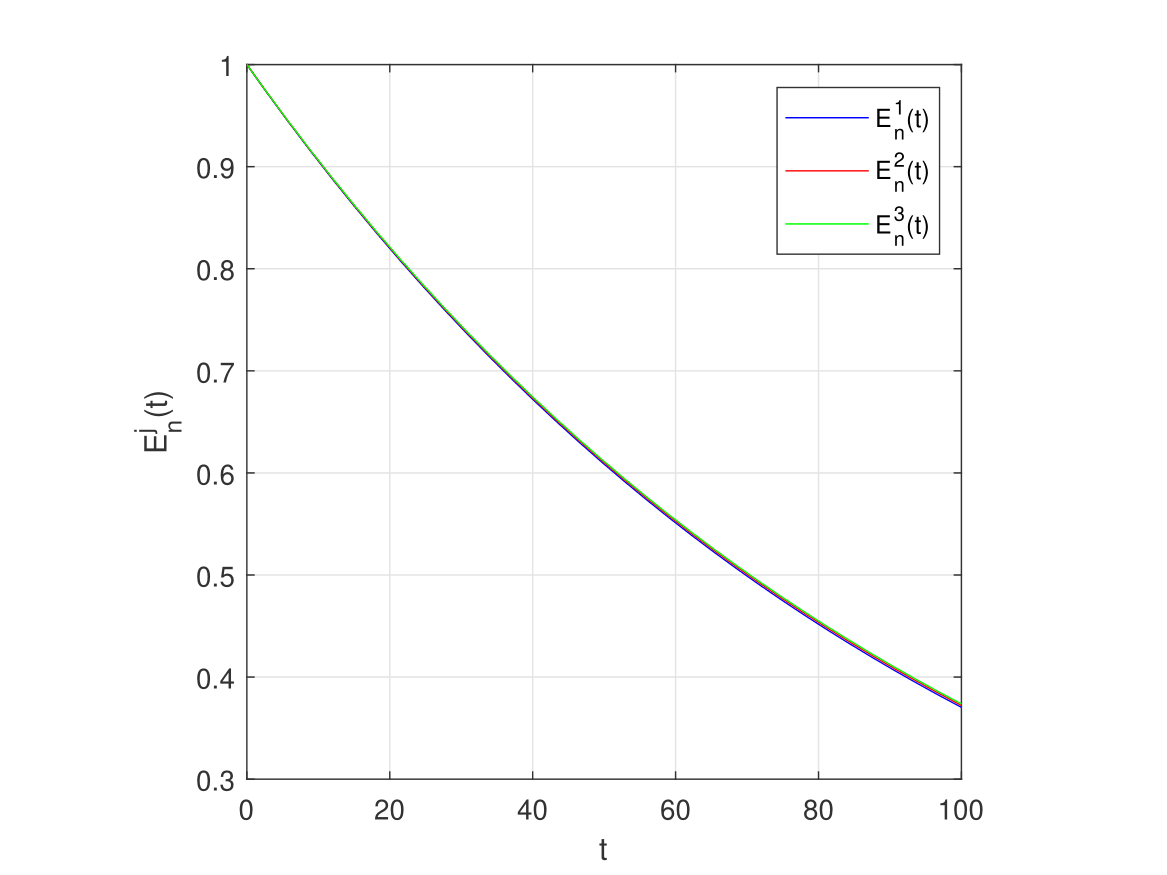}
\caption{No effect of smoothness on exponential decay of $E_{1,n}(t)$}
\end{center}
\end{figure}
\begin{figure}[H]
\begin{center}
\includegraphics[scale=0.8]{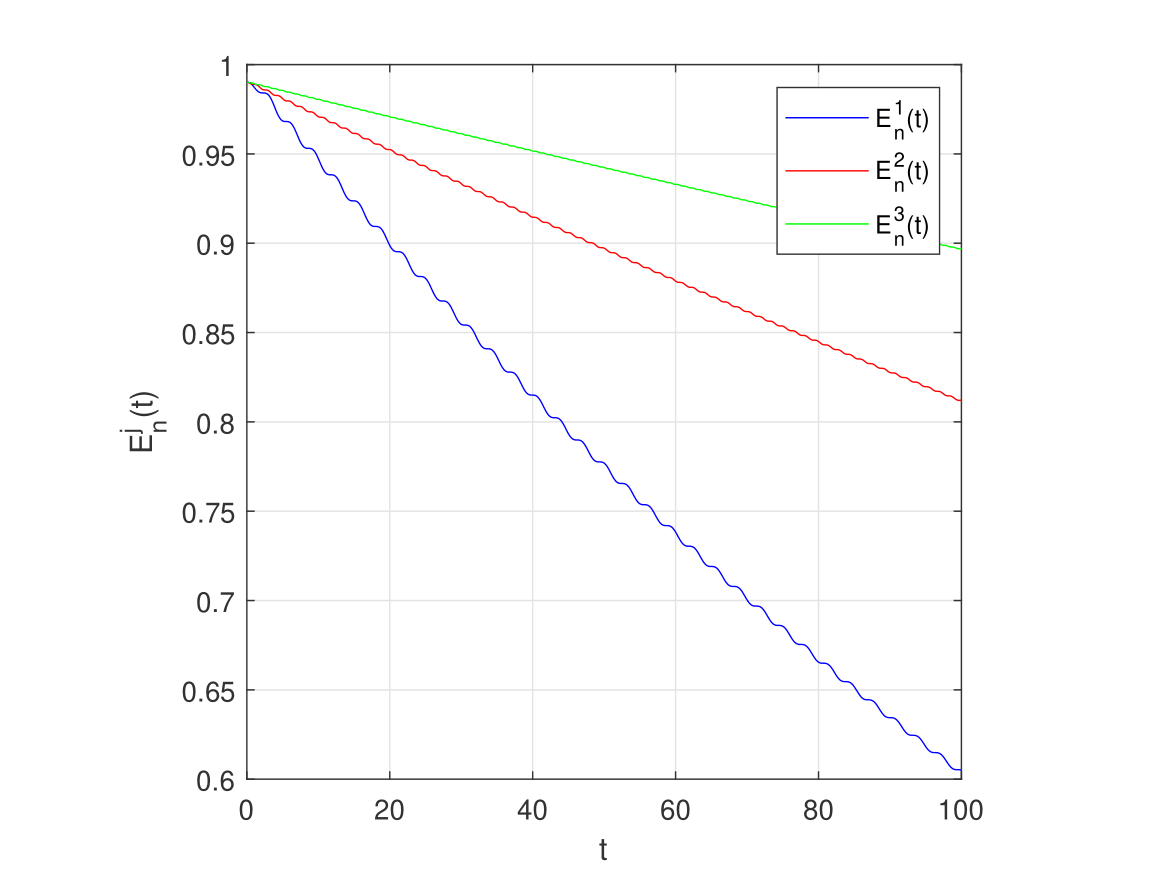}
\caption{Effect of smoothness on polynomial decay of $E_{2,n}(t)$}
\end{center}
\end{figure}

\section{Conclusions}
In this paper, we have investigated the question of homogenization and approximation of
thermoelasticity by using variational methods, semigroup approach and some space semi-discretization schemes. Once the appropriate setting is definied, we have shown the wellposedness of our system and we have recalled the notion of two scale convergence method as well as a weak version of Trotter Kato convergence Theorem.

To derive the homogenized thermoelastic system that captures the effective properties of the heterogeneous mediumn we inspire from \cite{KS2018} and use variational techniques and some regularity arguments. We have shown further that the solutions of our system decay polynomially and not exponentially to zero.

Afterwards, we have considered the approximation of two thermoelastic wave model with different asysmptotic behavior. One with exponential decay and the other one with polynomial decay. Due to numerical experiments, we have deduced that the two discrete systems show different behavior of the spectra and their discrete energies inherit the same behavior of the continuous ones. 

Finally we have shown numerically how the smoothness of data could impact the rate of decay of the energy associated the weakly coupled thermoelastic wave model, and we have concluded that these results are in perfect agreement with the theory already established in the  works \cite{LZb1994,MN16,N21}.

\section*{\bf Acknowledgements}

The author would like to thank Dr Shane Cooper for introducing him to the homogenization theory and for several helpful discussions, which, in turn, motivated him to consider the problem of this paper.

The author warmly thank the referee for a careful reading and useful comments, which improved the manuscript at several points.

\newpage

\section*{\bf Data availability statement}
No datasets were generated or analyzed during the current study.


\begin{thebibliography}{00}

\bibitem{Allaire92} G. Allaire, \textit{Homogenization and two-scale convergence}, SIAM J. Math. Anal., \textbf{23} (1992), 1482-1512.

\bibitem{Allaire93} G. Allaire, \textit{Two-Scale Convergence and Homogenization of Periodic Structures}, School on homogenization ICTP, Trieste, 1993.




\bibitem{BEPS2006} A. B\' atkai, K.J. Engel, J. Pr\"uss and R. Schnaubelt, \textit{Polynomial stability of operator semigroups}, Math. Nachr. 279 (2006) pp. 1425–1440.

\bibitem{BLP78} A. Bensoussan, J. L. Lions, G. Papanicolaou, \textit{Asymptotic Analysis for Periodic Structures}, North Holland, Amsterdam, 1978.

\bibitem{BT2010} A. Borichev and Y. Tomilov, \textit{Optimal polynomial decay of functions and operator semigroups}, Math. Ann., 347(2), pp.455-478,2010.

\bibitem{Carlson72} D. E. Carlson, \textit{Linear thermoelasticity}, in Handbuch der Physik, Ed.: C. Truesdell, Springer
Verlag, Berlin, 1972.

\bibitem{Chadwick60} P. Chadwick, \textit{Thermoelasticity. The dynamical theory, in Progress in Solid Mechanics}, vol. 1, Amsterdam, North-Holland, 1960.

\bibitem{CD00} D. Cioranescu, P. Donato, \textit{An Introduction to Homogenization}, Oxford University Press, Oxford, 2000.



\bibitem{EN2000} K. Engel and R. Nagel, \textit{One-parameter semigroups for linear evolution equations}. Encyclopedia of Mathematics and its Applications. Springer-Verlag, New York, 2000.

\bibitem{Evans88} L.C. Evans, \textit{Weak Convergence Methods for Nonlinear Partial Differential Equations}. AMS Regional
Conference Series in Mathematics, \textbf{74}, 1988.

\bibitem{Francfort83} G. A. Francfort, \textit{Homogenization and Linear Thermoelasticity}, SIAM Journal on Mathematical Analysis, SIAM J. Math. Anal., \textbf{14}, 1983, 696-708.

\bibitem{GRT92} J. S. Gibson, I. G. Rosen, and G. Tao, \textit{Approximation in control of thermoelastic systems}, SIAM J. Control
Optim., 30 (1992), pp. 1163-1189.

\bibitem{Hansen1992} S. W. Hansen, \textit{Exponential energy decay in a linear thermoelastic rod}. J. Math. Anal. Appli.,\textbf{167}, 1992, pp. 429-442.


\bibitem{JKO94} V. V. Jikov, S. M. Kozlov, and O. A. Oleinik, \textit{Homogenization of Differential Operators and Integral Functionals}, Springer, Berlin (1994).

\bibitem{KS2018} I. V. Kamotski $\&$ V. P. Smyshlyaev,\textit{Two-scale homogenization for a
general class of high contrast PDE systems with periodic coefficients}, Applicable Analysis, 2018. Doi: 10.1080/00036811.2018.1441994

\bibitem{KBT1997} F. A. Khodja, A. Benabdallah and D. Teniou, \textit{Dynamical stabilizers and coupled systems}, ESAIM Proceeding, \textbf{2}, (1997), pp. 253-262.

\bibitem{KBB1999} F. A. Khodja, A. Bader and A. Benabdallah, \textit{Dynamic stabilization of systems via decoupling techniques}, ESAIM: COCV, \textbf{4}, (1999), 577–593.

\bibitem{LR05} Z. Liu and B. Rao, \textit{Characterization of polynomial decay rate for the solution of linear evolution equation}. Z. Angew. Math. U. Phys. ZAMP 56 (2005) 630–644.

\bibitem{LZb1993} Z. Liu and S. Zheng, \textit{Exponential stability of semigroup associated with thermoelastic system}, Quart. Appl. Math, 51, (1993), pp. 535-545.

\bibitem{LZb1994} Z. Y. Liu and S. Zheng, \textit{Uniform exponential stability and approximation in control of a thermoelastic system}, SIAM J. Control Optim. 32, (1994), pp.1226-1246.

\bibitem{LZ99} Z. Y. Liu and S. Zheng, \textit{Semigroups Associated with Dissipative Systems}. Chapman \& Hall/CRC Research Notes in Mathematics Series (1999).




\bibitem{MN16} L. Maniar and S. Nafiri, \textit{Approximation and uniform polynomial stability of $C_{0}$ semigroups}, ESAIM: COCV, \textbf{22}, 208-235, 2016.

\bibitem{Milton2001} Graeme W. Milton, The Theory of Composites, Cambridge University Press (2001).

\bibitem{N21} S. Nafiri, \textit{Uniform Polynomial Decay and Approximation in Control of a Family of Abstract Thermoelastic Models}, JDCS, 2021.


\bibitem{Nechvatal03} L. Nechvátal, \textit{On two-scale convergence}, Mathematics and Computers in Simulation, \textbf{61}, 2003, 489–495.

\bibitem{Ng89} G. Nguetseng, \textit{A general convergence result for a functional related to the theory of homogenization}, SIAM J. Math. Anal. \textbf{20} (1989), 608-623. 

\bibitem{Ng90} G. Nguetseng, \textit{Asymptotic analysis for a stiff variational problem arising in mechanics}. SIAM J. Math. Anal. 21(6), 1394–1414 (1990)



\bibitem{P05} S. E. Pastukhova, \textit{On the Convergence of Hyperbolic Semigroups in Variable Hilbert Spaces}, Journal of Mathematical Sciences, \textbf{127}, 2005, 2263-2283.

\bibitem{Pazy83} A. Pazy, \textit{Semigroups of Linear Operators and Applications to Partial Differential Equations}, Springer-Verlag, New York, 1983.

\bibitem{Sanchez80} E. Sanchez-Palencia, \textit{Non-homogenous media and vibration theory}, Lecture Notes in Physics, \textbf{127}, Springer, Berlin, 1980.






\bibitem{Tartar78} L. Tartar, \textit{Topics in Nonlinear Analysis}, Publications mathématiques d’Orsay 78.13, Université de Paris-Sud, 1978.


\bibitem{Temam79} R. Temam, \textit{Navier-Stokes Equations}, North-Holland, 1979.



\bibitem{Yosida80} K. Yosida, \textit{Functional Analysis}, 6th ed., Springer-Vedag, New York, 1980.

\bibitem{ZP04}V. V. Zhikov, \textit{On Two-Scale Convergence}, Journal of Mathematical Sciences, \textbf{120}, 2004, 1328–1352.


\bibitem{ZP07} V. V. Zhikov and S. E. Pastukhova, \textit{On the Trotter–Kato Theorem in a Variable Space}, Functional Analysis and Its Applications, \textbf{41}, 2007, 264–270.

\end{thebibliography}
\end{document}